# Lie-Poisson groups and the Miura transformation

*JM Figueroa-O'Farrill* [1] [♮] *and S Stanciu* [2] [♯]

[♮] Queen Mary and Westfield College, London, UK
[♯] ICTP, Trieste, Italy

We point out that the recent proof of the Kupershmidt–Wilson theorem by Cheng and Mas–Ramos is underpinned by the Lie-Poisson property of the second Gel'fand–Dickey bracket. The supersymmetric Kupershmidt–Wilson theorem is also proved along these same lines. Finally we comment on the possible repercussions in the problem of the coproduct for W-algebras.

## Introduction

The relation between W-algebras and integrable hierarchies has been a close one almost since their inception. Essentially to every W-algebra there corresponds (via a sort of classical limit) a Poisson algebra of functions on the phase space of an integrable hierarchy of soliton equations. This correspondence has been exploited in order to define classical W-algebras, to quantise them, and to analyse their representation theory. In particular, one result in the theory of integrable hierarchies has played a prominent role in the applications of this correspondence: the Miura transformation [1] and the related Kupershmidt–Wilson theorem [2]. In a nutshell, the Kupershmidt–Wilson theorem underlies the free field realisations of the W-algebras which, together with the quantum Drinfel'd–Sokolov reduction [3], are the only known uniform constructions for quantum W-algebras. These free-field realisations underlie most of what is known about the representation theory of W-algebras.

Let us briefly review the content of the Kupershmidt–Wilson theorem. The notation and concepts about the Lax formalism and pseudodifferential operators (ΨDO's) are as in Dickey's book [4]. Recall that if $M_n$ denotes the affine space of Lax operators of the form $L = \partial^n + \sum_{i=1}^{n} u_i \partial^{n-i}$ then one can define a Poisson bracket on the differential ring $R[u]$ generated by the $\{u_i\}$:

$$\{u_i(x), u_j(y)\} = J_{ij}(x) \cdot \delta(x-y) , \qquad (1)$$

where the $J_{ij}(x)$ are differential operators whose coefficients lie in $R[u]$ which are given by the following construction. Let $L$ be as above, and we let $X = \sum_{i=1}^{n} \partial^{i-n-1} x_i$ be a ΨDO. Then one defines the $J_{ij}(x)$ by

$$J_L(X) = (LX)_+ L - L(XL)_+ = L(XL)_- - (LX)_- L \qquad (2)$$
$$= \sum_{i,j=1}^{n} (J_{ij} \cdot x_j) \partial^{n-i} .$$

The map $X \mapsto J_L(X)$ given by (2) is known as the Adler map and was introduced in [5] where it was conjectured that the bracket (1) was Poisson. This was proven by Gel'fand and Dickey [6] by direct calculation (see [4] for a cleaner version of the original proof) which led to their being called the (second) Gel'fand–Dickey brackets. In terms of W-algebras, they define the classical W-algebra known as $\mathsf{GD}_n$ out of which the more famous $\mathsf{W}_n$ algebra is obtained by hamiltonian reduction according to the constraint $u_1 = 0$.

At least three more proofs are known of the Poisson property of (1). On the one hand, Drinfel'd and Sokolov [3], based on previous work of Reiman and Semenov-Tyan-Shanskiĭ, proved that (1) arises via hamiltonian reduction from the natural Poisson structure in the dual of the loop algebra of $GL(n)$. On the other hand, and more to the point of this letter, the Kupershmidt–Wilson theorem [2] states that (1) is induced from the Poisson bracket in the symmetric algebra of the $n$-dimensional Heisenberg algebra via the Miura transformation. Explicitly, suppose that we formally factorise the Lax operator $L$ as follows:

$$L = (\partial + \phi_1)(\partial + \phi_2) \cdots (\partial + \phi_n) . \qquad (3)$$

This factorisation induces an embedding of differential rings $\mu : R[u] \to R[\phi]$ which generalises the Miura transformation [1]. If on $R[\phi]$ we define the following Poisson structure,

$$\{\phi_i(x), \phi_j(y)\} = \delta_{ij} \delta'(x-y) , \qquad (4)$$

then we can compute the induced bracket on the image of the Miura transformation $u_i(\phi)$. The Kupershmidt–Wilson theorem states that the induced bracket is none other than the second Gel'fand–Dickey bracket (1). This theorem hides within it two equally remarkable facts: first that the induced brackets of the $u_i(\phi)$ should once again close on the $u_i(\phi)$, and secondly that it should precisely

[1] e-mail: j.m.figueroa@qmw.ac.uk
[2] e-mail: sonia@ictp.trieste.it

agree with (1). Finally, on the third hand, Semenov-Tyan-Shanskiĭ [7] proved that the second Gel'fand–Dickey bracket is the Drinfel'd–Sklyanin bracket associated with a classical $R$-matrix in the Lie algebra $\mathbb{A}$ of $\Psi$DO's. If we split $\mathbb{A} = \mathbb{A}_+ \oplus \mathbb{A}_-$ into the subalgebras of differential and Volterra operators respectively, and we let $P_\pm : \mathbb{A} \to \mathbb{A}$ denote the projectors onto $\mathbb{A}_\pm$ along $\mathbb{A}_\mp$, then the relevant $R$-matrix $R : \mathbb{A} \to \mathbb{A}$ is given by $R = P_+ - P_-$. The Drinfel'd–Sklyanin bracket associated with an $R$-matrix, is given (for the case of an associative algebra $\mathbb{A}$ with a nondegenerate trace Tr) by

$$\{F,G\}(L) = \operatorname{Tr} dFLR(dGL) - \operatorname{Tr} LdFR(LdG) . \tag{5}$$

It is easy to see that the Gel'fand–Dickey bracket is recovered (up to an overall factor) if we choose $R = P_+ - P_-$. In this paper we will see that in fact these two proofs of the Poisson property of the Gel'fand–Dickey bracket are intimately related.

Several proofs of the Kupershmidt–Wilson theorem exist in the literature, each emphasising one or both of the two remarkable facts above. Kupershmidt's and Wilson's original proof [2] was both fairly long and not elementary. A much shorter yet elementary proof was subsequently given by Dickey [8], using only algebraic manipulations with $\Psi$DO's. A third, more indirect proof is provided by the Drinfel'd–Sokolov reduction scheme [3]. Not unrelated to this approach Wilson [9] presented a fourth proof—or rather a conceptualisation which addresses the first remarkable fact—of the Kupershmidt–Wilson theorem in terms of differential Galois theory and introducing in the process an integrable hierarchy underlying both the KdV- and modified KdV-type hierarchies. Finally, a fifth proof has appeared recently in the independent work of Y Cheng [10] and Mas–Ramos [11], which exploits the behaviour of the Adler map (2) under multiplication of Lax operators.

The point of the present letter is to conceptualise the results of [10] and [11] in terms of (a slight generalisation of) Drinfel'd's theory of Lie-Poisson groups. Explicitly, we show that underpinning the results of Cheng and Mas–Ramos, lies the fact that the space $M = \cup_n M_n$ is a (formal) Lie-Poisson group. This fact itself is not new, having appeared originally in Semenov-Tyan-Shanskiĭ's work [7] on classical $R$-matrices. What is novel in this letter, however, is its relation to the Kupershmidt–Wilson theorem. In fact, it seems that the Kupershmidt–Wilson theorem typifies and epitomises the Lie-Poisson property.

Before entering into the details, let us comment on an amusing historical annecdote. In the introduction of the paper of Kupershmidt and Wilson [2] they state that, as a result of their theorem, the second Gel'fand–Dickey bracket is characterized by two properties: that for a linear differential operator $L = \partial + \phi_i$ it is simply given by $\partial$, and that multiplication of operators is a "canonical transformation." At that time, of course, Lie-Poisson theory had not been formulated and as a result their proof of their theorem was more involved than it might have otherwise been.

Lie-Poisson basics

We start by reviewing the basic concepts in the theory of Lie-Poisson groups [12] needed to put this work in the right context. A Lie-Poisson group is a Lie group $G$ which is also a Poisson manifold and such that the group operations and the Poisson bracket are compatible. The relevant compatibility condition can be stated in many ways. At the algebraic level, one demands that the algebra of functions $\operatorname{Fun}(G)$ on the group be a Poisson-Hopf algebra. That is, one demands that the comultiplication $\Delta : \operatorname{Fun}(G) \to \operatorname{Fun}(G) \otimes \operatorname{Fun}(G)$ be a Poisson morphism. The algebraic approach is not very fruitful in our case since the formal geometry with which we endow the space of Lax operators does not afford the functions with an associative structure but only with a Poisson bracket. We will find it more convenient to argue dually and rephrase the Lie-Poisson condition. Under the isomorphism $\operatorname{Fun}(G) \otimes \operatorname{Fun}(G) \cong \operatorname{Fun}(G \times G)$, the comultiplication is induced by the group multiplication $m : G \times G \to G$. The group is Lie-Poisson provided that group multiplication is a Poisson morphism, where we endow $G \times G$ with the product Poisson structure. In other words, $G$ is Lie-Poisson provided that for all functions $f, g \in \operatorname{Fun}(G)$, and all $x, y \in G$,

$$\{m^* f, m^* g\}(x,y) = m^* \{f,g\}(x,y) = \{f,g\}(xy) . \tag{6}$$

If the Poisson bracket satisfies this identity it is said to be *grouped* or *multiplicative*. We prefer this latter name since we will not be always working with groups.

To describe this more explicitly, it is necessary to look more closely at the product Poisson structure in $G \times G$. Let $\operatorname{pr}_i : G \times G \to G$, for $i = 1, 2$, be the Cartesian projection on the $i$th factor. The product Poisson structure on $G \times G$ is uniquely characterised by the following two properties: (1) $\operatorname{pr}_i$ are Poisson morphisms; and (2) $\{\operatorname{pr}_1^* f, \operatorname{pr}_2^* g\} = 0$ for all functions $f, g \in \operatorname{Fun}(G)$. One can verify that if $\bar{f}, \bar{g}$ are functions in $G \times G$ then their Poisson bracket is given by

$$\{\bar{f}, \bar{g}\}(x,y) = \{\ell_x^* \bar{f}, \ell_x^* \bar{g}\}(y) + \{r_y^* \bar{f}, r_y^* \bar{g}\}(x) , \tag{7}$$

where $r_y : G \to G \times G$ and $\ell_x : G \to G \times G$ are defined by $r_y(z) = (z,y)$ and $\ell_x(z) = (x,z)$. Using (7) into (6), we find that a Poisson bracket on $G$ is multiplicative if and only if it satisfies

$$\{f,g\}(xy) = \{\lambda_x^* f, \lambda_x^* g\}(y) + \{\varrho_y^* f, \varrho_y^* g\}(x) , \tag{8}$$

where we $\lambda_x = m \circ \ell_x : G \to G$ and $\varrho_y = m \circ r_y : G \to G$ are left and right multiplications by $x$ and $y$ respectively.

Notice that equation (8) makes sense even if $G$ is not a group but only a groupoid. In other words, it is not necessary that there should exist a unit for the multiplication nor an inverse. In summary, we can consider Poisson manifolds $M$ with an associative multiplication $M \times M \to M$ with unit[3]; and a Poisson structure on $M$ is said to be multiplicative if condition (8) is satisfied. The invertible elements of $M$ will then form a Lie-Poisson group. This represents a mild generalisation of Drinfel'd's theory of Lie-Poisson groups which has recently been advocated in [13].

### Lie-Poisson property of the Adler map

We now apply these considerations to the space of Lax operators. Let $F, G$ be functions on $M$ and let us consider the class of Poisson brackets which can be expressed in the form

$$\{F, G\}(L) = \operatorname{Tr} H_L(dF) dG \ , \tag{9}$$

for some hamiltonian map $H_L$ and where Tr denotes the (essentially unique) trace on the space of $\Psi$DO's introduced by Adler in [5]. According to the previous section, the Poisson bracket (9) will be multiplicative provided that the following identity is satisfied:

$$\{F, G\}(AB) = \{\varrho_B^* F, \varrho_B^* G\}(A) + \{\lambda_A^* F, \lambda_A^* G\}(B) \ . \tag{10}$$

This identity will in turn be satisfied provided that the hamiltonian map $H_L$ obeys an identity relating $H_{AB}$ to $H_A$ and $H_B$. First of all notice that for any $\Psi$DO $V$,

$$\operatorname{Tr} d(\lambda_A^* F) V = \operatorname{Tr}(\lambda_A^* dF) V = \operatorname{Tr} dF (\lambda_A)_* V = \operatorname{Tr} dF A V \ , \tag{11}$$

whence $d\lambda_A^* F = dF A$. Similarly, $d\varrho_B^* F = B dF$. Putting this all together into (10), we find that (10) is satisfied for all $F$ and $G$ if and only if

$$H_{AB}(X) = H_A(BX) B + A H_B(XA) \ , \tag{12}$$

for all $\Psi$DO's $A, B, X$; or more intrinsically,

$$H_{AB} = (\varrho_B)_* \circ H_A \circ \varrho_B^* + (\lambda_A)_* \circ H_B \circ \lambda_A^* \ . \tag{13}$$

It is now a simple calculation to show that the Adler map does indeed obey this identity:

$$J_{AB}(X) = (ABX)_+ AB - AB(XAB)_+$$

---

[3] the unit is not essential, but it exists in the examples we will be interested in.

$$= (ABX)_+ AB \pm A(BXA)_+ B - AB(XAB)_+$$
$$= A\left[(BXA)_+ B - B(XAB)_+\right] + \left[(ABX)_+ A - A(BXA)_+\right] B$$
$$= A J_B(XA) + J_A(BX) B \ . \tag{14}$$

In other words, the second Gel'fand–Dickey bracket is multiplicative. As we now see this property underpins the proof in [11] (see also [10]) of the Kupershmidt–Wilson theorem.

It is easier in fact, to go back to the original definition of a multiplicative Poisson bracket, as one which under which the multiplication is a Poisson morphism. Let $L = \partial^n + \sum_{i=1}^n u_i \partial^{n-i} = AB$ where $A = \partial^p + \sum_{i=1}^p a_i \partial^{p-i}$ and $B = \partial^q + \sum_{i=1}^q b_i \partial^{q-i}$ and $p + q = n$. A function $F(L) = \int f(u)$ can be considered as a function $F(A, B) = \int f(a, b)$. (We abuse the notation slightly by not changing the name of the function. In fact, the two-variable function $F$ is nothing but the pull-back via the multiplication of the one-variable function $F$.) The multiplicative property of the second Gel'fand–Dickey bracket is then equivalent to

$$\{F, G\}(L) = \sum_{i,j=1}^n \frac{\delta F}{\delta u_i} \cdot J_{ij}^{(n)} \cdot \frac{\delta G}{\delta u_j}$$
$$= \sum_{i,j=1}^p \frac{\delta F}{\delta a_i} \cdot J_{ij}^{(p)} \cdot \frac{\delta G}{\delta a_j} + \sum_{i,j=1}^q \frac{\delta F}{\delta b_i} \cdot J_{ij}^{(q)} \cdot \frac{\delta G}{\delta b_j} \ . \tag{15}$$

Iterating this theorem to the complete factorisation $L = A_1 A_2 \cdots A_n$ of $L$ into linear factors $A_i = \partial + \phi_i$, yields

$$\sum_{i,j=1}^n \frac{\delta F}{\delta u_i} \cdot J_{ij}^{(n)} \cdot \frac{\delta G}{\delta u_j} = \sum_{i=1}^n \frac{\delta F}{\delta \phi_i} \cdot J^{(1)} \cdot \frac{\delta G}{\delta \phi_i} \ , \tag{16}$$

which is the Kupershmidt–Wilson theorem after noticing that $J^{(1)} = \partial$. Notice that at the level of W-algebras, (15) says that $\mathsf{GD}_n$ is embedded in $\mathsf{GD}_p \times \mathsf{GD}_q$ for $p + q = n$.

### Some generalisations

The above results admit several straightforward generalisations. First of all one can consider $q\Psi$DO's [14]. Results along these lines have been obtained by Cheng [10] and we shall have nothing further to add to this here. Another line of generalisation is to consider supersymmetric $\Psi$DO's [15] [16]. This has the virtue of yielding a more constructive and conceptual proof of the supersymmetric Kupershmidt–Wilson theorem first proven in [17]. We do this now briefly, following the conventions and notation of [16].

Supersymmetric $\Psi$DO's take the form $L = \sum_i U_i D^i$ where the $U_i$ are superfields in a $(1|1)$ superspace $U_i = \psi_i + \theta u_i$, where $\theta^2 = 0$, and where $D$ denotes the supercovariant derivative $D = \frac{\partial}{\partial \theta} + \theta \partial$, which obeys $D^2 = \partial$. We take $L$ to be homogeneous of degree $|L|$ under the $\mathbb{Z}_2$-grading. We shall be considering functions of the form $F[U] = \int_B f(U)$ with $\int_B$ the Berezin integral where $f$ is also homogeneous under the $\mathbb{Z}_2$ grading. The supersymmetric Gel'fand–Dickey bracket of two homogeneous functions is given by

$$\{F,G\}(L) = -(-)^{|F|+|G|+|L|} \mathrm{Str} J_L(d_L F) d_L G , \qquad (17)$$

where the Adler map $J_L$ is given by the same expression as in (2) with the obvious reinterpretation of the symbols and we write $d_L F$ for the gradient of $F$ as a function of $L$. Under a vector field $\delta L$, a function changes by

$$\delta F = -(-)^{|\delta L|+|F|} \mathrm{Str} \delta L dF . \qquad (18)$$

If $L = AB$, then

$$\delta L = \delta A B + (-)^{|\delta||A|} A \delta B , \qquad (19)$$

whence if we write $d_A F$ and $d_B F$ for the gradients of $F[AB]$ as functions of $A$ and $B$ respectively, one finds

$$\begin{aligned} d_A F &= (-)^{|A|+|L|} B d_L F \\ d_B F &= (-)^{|A||F|} d_L F A . \end{aligned} \qquad (20)$$

Plugging this into the expression of the supersymmetric Gel'fand–Dickey bracket (17), one finds for $L = AB$,

$$\{F,G\}(L) = -(-)^{|F|+|G|+|L|} \left( \mathrm{Str} J_A(d_A F) d_A G + (-)^{|A|} \mathrm{Str} J_B(d_B F) d_B G \right) . \qquad (21)$$

The only remarkable point here is the relative sign between the two terms which explains the choice of starting Poisson bracket for the $\Phi_i$ in the Miura transformation $L = D^n + \sum_{i=1}^n U_i D^{n-i} = (D+\Phi_1)(D+\Phi_2)\cdots(D+\Phi_n)$, which was used in [17] to formulate and prove the supersymmetric Kupershmidt–Wilson theorem.

Finally, let us mention that the results on the Lie-Poisson structure generalise to any associative algebra $\mathbb{A}$ admitting a nondegenerate trace Tr and a vector space decomposition into two subalgebras $\mathbb{A} = \mathbb{A}_+ \oplus \mathbb{A}_-$, with the additional property that $\mathbb{A}_\pm$ are maximally isotropic under the bilinear form defined by the trace $\langle a , b \rangle = \mathrm{Tr}\, ab$. Indeed, under the commutator, $\mathbb{A}$ becomes a Lie algebra and $\mathbb{A}_\pm$ are Lie subalgebras; furthermore the bilinear form $\langle - , - \rangle$ is ad-invariant. In other words, therefore $(\mathbb{A}, \mathbb{A}_+, \mathbb{A}_-)$ is a Manin triple and $\mathbb{A}$ and $\mathbb{A}_\pm$ are Lie bialgebras. In particular, $\mathbb{A}$ is triangular, with $R$-matrix given by $R = P_+ - P_-$; where $P_\pm$ are the projectors onto $\mathbb{A}_\pm$ along $\mathbb{A}_\mp$. This implies that the (formal) Lie groups associated with $\mathbb{A}_\pm$ and $\mathbb{A}$ are Lie-Poisson relative to the Poisson bracket associated with $R$ (see [12]). For the particular case of $\mathbb{A}$ being the algebra of $\Psi$DO's, Semenov-Tyan-Shanskiĭ showed that the Poisson bracket associated with $R$ agrees precisely with the second Gel'fand–Dickey bracket (1).

## Outlook

We end the paper with two remarks: a remark on the uniqueness of the Adler map and a remark on the possible implications of this result for the representation theory of W-algebras.

For all its remarkable properties, the Adler map (2) still remains shrouded in the aura of mystery usually associated with constructions which owe more to inspiration than to a well-axiomatised formalism. One could hope that the fact that it is multiplicative relative to the multiplication of Lax operators would go some way towards this aim; but this is of course not the case. If we restrict ourselves to the space $M = \cup_n M_n$—consisting of Lax operators which are formally invertible in the ring of $\Psi$DO's and hence which define a formal Lie group, we can apply standard Drinfel'd theory to associate to every multiplicative Poisson bracket, a 1-cocycle in the Lie algebra $\mathbb{A}$ of $\Psi$DO's with values in linear maps $\mathbb{A} \to \mathbb{A}$. Assuming that this 1-cocycle is a coboundary by an $R$-"matrix" $R : \mathbb{A} \to \mathbb{A}$, then formula (5) defines a Poisson bracket which share with the second Gel'fand–Dickey bracket the same hamiltonians. It would be interesting to find other $R$-matrices.

Finally, despite the progress brought about by the applications of the techniques of the Lax formalism to the study of W-algebras, one aspect remains singularly unsolved: the existence of tensor product representations. Given two linear representations of a Lie algebra, the tensor product of their underlying vector spaces also admits a representation inherited naturally from the representations of each of the tensorands. Underlying this property lies the algebraic fact that the universal envelope of a Lie algebra is a Hopf algebra and, in particular, possesses a coproduct. Unlike Lie algebras, the intrinsic nonlinearity present in a W-algebra forbids the construction of the tensor product representation in the same way as for a Lie algebra, and indeed to this day, there is no known coproduct for W-algebras. From the interpretation of classical W-algebras as algebras of functions on the space of Lax operators, it would make sense to try and glean some information about the existence of a coproduct in this setting using the multiplication present in the space of Lax operators; after all, the coproduct in the universal envelope of a Lie algebra is induced by the multiplication on the Lie group. From what we see above, it seems that GD$_n$ does not inherit a coproduct from the multiplication of the Lax operators, since this represents a map GD$_n \to$ GD$_p \times$ GD$_q$ ($p+q=n$). It seems that one must work with all the Gel'fand–Dickey algebras at the same time. It is not clear to us at the present stage how to implement this in terms of conformal field theory.

## ACKNOWLEDGEMENTS

We are grateful to Javier Mas for communicating the results of [11] prior to publication, and to Yi Cheng for sending us his preprint [10]. In addition, SS would like to thank ICTP (Trieste) for the leave of absence to visit QMW College; and Chris Hull and the String Theory Group at QMW for the hospitality extended to her during the completion of this work.